\documentclass[titlepage,12pt]{article}

\usepackage{graphicx}
\usepackage{amssymb}
\usepackage{psfrag}

\newcommand{\R}{\mathbb R}

\begin{document}

\baselineskip=18pt

\begin{center}
{\Large{\bf Hopf Bifurcations in a
 \\ Watt Governor with a Spring}}
\end{center}

\begin{center}
{\large Jorge Sotomayor}
\end{center}
\begin{center}
{\em Instituto de Matem\'atica e Estat\'{\i}stica, Universidade de
S\~ao Paulo\\ Rua do Mat\~ao 1010, Cidade Universit\'aria\\ CEP
05.508-090, S\~ao Paulo, SP, Brazil
\\}e--mail:sotp@ime.usp.br
\end{center}
\begin{center}
{\large Luis Fernando Mello}
\end{center}
\begin{center}
{\em Instituto de Ci\^encias Exatas, Universidade Federal de
Itajub\'a\\Avenida BPS 1303, Pinheirinho, CEP 37.500-903,
Itajub\'a, MG, Brazil
\\}e--mail:lfmelo@unifei.edu.br
\end{center}
\begin{center}
{\large Denis de Carvalho Braga}
\end{center}
\begin{center}
{\em Instituto de Sistemas El\'etricos e Energia, Universidade
Federal de Itajub\'a\\Avenida BPS 1303, Pinheirinho, CEP
37.500-903, Itajub\'a, MG, Brazil
\\}e--mail:braga@unifei.edu.br
\end{center}

\begin{center}
{\bf Abstract}
\end{center}

\vspace{0.1cm}

This paper pursues the study carried out by the authors in {\it
Stability and Hopf bifurcation in a hexagonal governor system}
\cite{smb3}, focusing on the codimension one Hopf bifurcations in
the hexagonal Watt governor differential system. Here are studied
the codimension two, three and four Hopf bifurcations and the
pertinent Lyapunov stability coefficients and bifurcation
diagrams, illustrating the number, types and positions of
bifurcating small amplitude periodic orbits, are determined. As a
consequence it is found an open region in the parameter space
where two attracting periodic orbits coexist with an attracting
equilibrium point.

\vspace{0.1cm}



\newpage
\section{\bf Introduction}\label{intro}

The centrifugal governor is a device that automatically controls
the speed of an engine. The most important one is due to James
Watt --Watt governor-- and it can be taken as the starting point
for automatic control theory. Centrifugal governor design received
several important modifications as well as other types of
governors were also developed. From MacFarlane \cite{mac}, p. 251,
we quote:

\noindent ``Several important advances in automatic control
technology were made in the latter half of the 19th century. A key
modification to the flyball governor was the introduction of a
simple means of setting the desired running speed of the engine
being controlled by balancing the centrifugal force of the
flyballs against a spring, and using the preset spring tension to
set the running speed of the engine".

In this paper the system coupling the Watt governor with a spring
(resp. Watt governor) and the steam engine will be called simply
the Watt Governor System with Spring (WGSS) (resp. Watt Governor
System (WGS)). The stability analysis of the stationary states and
small amplitude oscillations of this  system will be pursued here.

The first mathematical analysis of the stability conditions in the
WGS was due to Maxwell \cite{max} and, in a user friendly style
likely to be better understood by engineers, by Vyshnegradskii
\cite{vysh}. A simplified version of the WGS local stability based
on the work of Vyshnegradskii is presented by Pontryagin
\cite{pon}.

From the mathematical point of view, the oscillatory, small
amplitude, behavior in the WGS can be associated to a periodic
orbit that appears from a Hopf bifurcation. This was established
by Hassard et al. in \cite{has1}, Al-Humadi and Kazarinoff in
\cite{humadi} and by the authors in \cite{smb1, smb2}. Another
procedure, based in the method of harmonic balance, has been
suggested by Denny \cite{denny} to detect large amplitude
oscillations.

In \cite{smb1} we characterized the surface of Hopf bifurcations
in a WGS, which is more general than that presented by Pontryagin
\cite{pon}, Al-Humadi and Kazarinoff \cite{humadi} and Denny
\cite{denny}.

In \cite{smb2} restricting ourselves to Pontryagin's system of
differential equations for the WGS, we carried out a deeper
investigation of the stability of the equilibrium along the critical
Hopf bifurcations up to codimension 3, happening at a unique point
at which the bifurcation diagram was established. A conclusion
derived from the diagram implied the existence of parameters where
the WGS has an attracting periodic orbit coexisting with an
attracting equilibrium.

In \cite{smb3} we characterized the hypersurface of Hopf
bifurcations in a WGSS. See Theorem \ref{teoremaHexa} and Fig.
\ref{sinalL1} for a review of the critical surface where the first
Lyapunov coefficient vanishes.

In the present paper we go deeper investigating the stability of the
equilibrium along the above mentioned critical surface. To this end
the second Lyapunov coefficient is calculated and it is established
that it vanishes along two curves. The third Lyapunov coefficient is
calculated on these curves and it is established that it vanishes at
a unique point. The fourth Lyapunov coefficient is calculated at
this point and found to be negative. See Theorem \ref{crucial}. The
pertinent bifurcation diagrams are established. See Fig. \ref{PointH} and
\ref{PointH1}. A conclusion derived from these diagrams, concerning the
region ---a solid ``tongue"--- in the space of parameters where two
attracting periodic orbits coexist with an attracting equilibrium,
is specifically commented in Section \ref{conclusion}.

The extensive calculations involved in Theorem \ref{crucial} have
been corroborated with the software MATHEMATICA 5 \cite{math} and
the main steps have been posted in the site \cite{mello}.

This paper is organized as follows. In Section \ref{S2} we introduce
the differential equations that model the WGSS. The stability of the
equilibrium point of this model is analyzed and a general version of
the stability condition is obtained and presented in the terminology
of Vyshnegradskii. The Hopf bifurcations in the WGSS differential
equations are studied in Sections \ref{S3} and \ref{S4}. Expressions
for the second, third and fourth Lyapunov coefficients, which fully
clarify their sign, are obtained, pushing forward the method found
in the works of Kuznetsov \cite{kuznet, kuznet2}. With this data,
the bifurcation diagrams are established. Concluding comments,
synthesizing and interpreting the results achieved here, are
presented in Section \ref{conclusion}.

\section{The Watt governor system with spring}\label{S2}

\newtheorem{teo}{Theorem}[section]
\newtheorem{lema}[teo]{Lemma}
\newtheorem{prop}[teo]{Proposition}
\newtheorem{cor}[teo]{Corollary}
\newtheorem{remark}[teo]{Remark}
\newtheorem{example}[teo]{Example}

\subsection{WGSS differential equations}\label{S21}

\begin{figure}[!h]
\centerline{
\includegraphics[width=12cm]{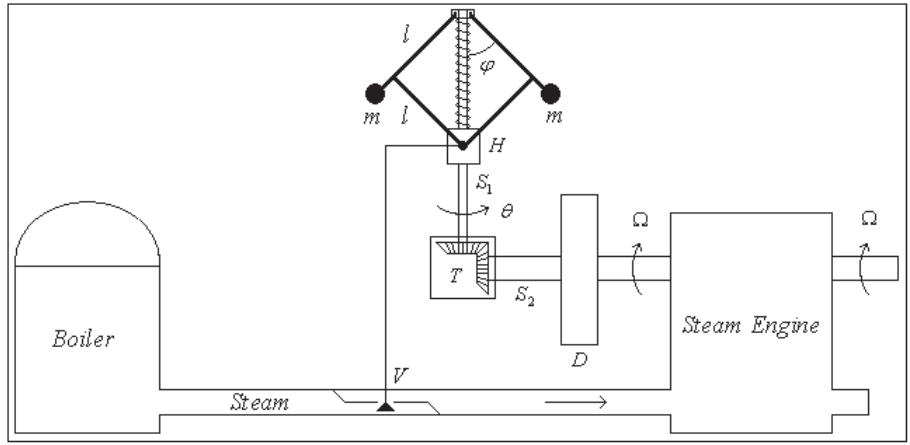}}
\caption{{\small Watt centrifugal governor with a spring -- steam
engine system}.}

\label{spring}
\end{figure}

The WGSS studied in this paper is shown in Fig. \ref{spring}.
There, $\varphi \in \left( 0,\frac{\pi}{2} \right)$ is the angle
of deviation of the arms of the governor from its vertical axis
$S_1$, $\Omega \in [0,\infty)$ is the angular velocity of the
rotation of the engine flywheel $D$, $\theta$ is the angular
velocity of the rotation of $S_1$, $l$ is the length of the arms,
$m$ is the mass of each ball, $H$ is a  sleeve which supports the
arms and slides along $S_1$, $T$ is a set of transmission gears
and $V$ is the valve that determines the supply of steam to the
engine.

The WGSS differential equations can be found as follows. For
simplicity, we neglect the mass of the sleeve and the arms. There
are four forces acting on the balls at all times. They are the
tangential component of the gravity
\[
-m g \sin \varphi,
\]
where $g$ is the standard acceleration of gravity; the tangential
component of the centrifugal force
\[
m \; l \; \sin \varphi \; \theta^2 \cos \varphi;
\]
the tangential component of the restoring force due to the spring
\[
-2 k l (1 - \cos \varphi) \sin \varphi,
\]
$2l$ is the natural length of the spring and $k \geq 0$ is the
spring constant; and the force of friction
\[
- b l \dot {\varphi},
\]
$b > 0$ is the friction coefficient.

From the Newton's Second Law of Motion, and using the transmission
function $\theta = c \: \Omega$, where $c > 0$, one has
\begin{equation}
\ddot{\varphi} = \left( \frac{2k}{m} + c^2 \Omega^2 \right)  \sin
\varphi \cos \varphi - \frac{2kl + mg}{ml} \: \sin \varphi
-\frac{b}{m} \: \dot{\varphi}. \label{Newton}
\end{equation}

The torque acting upon the flywheel $D$ is
\begin{equation}
I \: \dot{\Omega} = \mu \: \cos \varphi - F, \label{torque}
\end{equation}
where $I$ is the moment of inertia of the flywheel, $F$ is an
equivalent torque of the load and $\mu > 0$ is a proportionality
constant. See \cite{pon}, p. 217, for more details.

From Eq. (\ref{Newton}) and (\ref{torque}) the differential
equations of our model are given by
\begin{eqnarray}\label{sistema}
\frac{d \; \varphi}{d \tau} &=& \psi \nonumber\\
\frac{d \; \psi}{d \tau} &=& \left( \frac{2k}{m} + c^2 \Omega^2
\right)  \sin \varphi \cos \varphi -
\frac{2kl + mg}{ml} \; \sin \varphi -\frac{b}{m} \; \psi \\
\frac{d \; \Omega}{d \tau} &=& \frac{1}{I} \; \left( \mu \cos
\varphi -F \right) \nonumber
\end{eqnarray}
where $\tau$ is the time.

The standard Watt governor differential equations in Pontryagin
\cite{pon}, p. 217,
\begin{eqnarray}\label{standard}
\frac{d \; \varphi}{d \tau} &=& \psi \nonumber\\
\frac{d \; \psi}{d \tau} &=& c^2 \; \Omega^2 \; \sin \varphi \;
\cos \varphi -
\frac{g}{l} \; \sin \varphi - \frac{b}{m} \; \psi \\
\frac{d \; \Omega}{d \tau} &=& \frac{1}{I} \; \left( \mu \cos
\varphi -F \right) \nonumber
\end{eqnarray}
are obtained from (\ref{sistema}) by taking $k = 0$.

Defining the following changes in the coordinates, parameters and
time
\[
x = \varphi, \: y = \left( \frac{ml}{2kl+mg} \right)^{1/2} \psi,
\: z = c \left( \frac{ml}{2kl+mg} \right)^{1/2} \Omega,
\]
\[
t = \left( \frac{2kl+mg}{ml} \right)^{1/2} \tau, \: \kappa =
\frac{2kl}{2kl+mg},
\]
\[
\varepsilon = \frac{b}{m} \left( \frac{ml}{2kl+mg} \right)^{1/2},
\: \alpha = \frac{c \mu}{I} \left( \frac{ml}{2kl+mg} \right), \:
\beta = \frac{F}{\mu},
\]
where $0 \leq \kappa < 1$, $\varepsilon > 0$, $\alpha
> 0$ and $0 < \beta <1$, the differential equations
(\ref{sistema}) can be written as
\begin{eqnarray}\label{sistemafinal}
x' = \frac{d x}{d t} &=& y \nonumber\\
y' = \frac{d y}{d t} &=& (z^2 + \kappa) \; \sin x \; \cos x -
\sin x - \varepsilon \; y \\
z' = \frac{d z}{d t} &=& \alpha \; (\cos x - \beta) \nonumber
\end{eqnarray}
or equivalently by
\begin{equation}
{\bf x}' = f({\bf x}, {\bf \zeta}), \label{campo}
\end{equation}
where
\[
f({\bf x},{\bf \zeta}) = \left( y, (z^2 + \kappa) \; \sin x \;
\cos x - \sin x - \varepsilon \; y, \alpha \; \left( \cos x -
\beta \right) \right),
\]
\[
{\bf x} =(x,y,z) \in \left( 0, \frac{\pi}{2} \right) \times \R
\times [0,\infty)
\]
and
\[
{\bf \zeta} = (\beta, \alpha,\varepsilon, \kappa) \in \left( 0,
1\right) \times \left( 0, \infty \right) \times \left(0,\infty
\right) \times \left[ 0, 1 \right).
\]

\subsection{Stability analysis of the equilibrium
point}\label{S22}

The WGSS differential equations (\ref{sistemafinal}) have only one
admissible equilibrium point
\begin{equation}
P_0 = (x_0 , y_0 , z_0) = \left( \arccos \beta, 0, \left (
\frac{1}{\beta} - \kappa \right)^{1/2} \right). \label{P0}
\end{equation}

The Jacobian matrix of $f$ at $P_0$ has the form
\begin{equation}
Df \left( P_0 \right) = \left( \begin{array}{ccc}
0 & 1 & 0 \\
\\- \omega_0 ^2 & - \varepsilon & \xi \\
\\ -\alpha (1 - \beta^2)^{1/2}  & 0  & 0
\end{array} \right),
\label{jacobianP0}
\end{equation}
where
\begin{equation}
\omega_0 = \sqrt {\frac{1 - \beta ^2}{\beta}} \label{omega0}
\end{equation}
and
\[
\xi = 2 \beta^{1/2} (1 - \beta ^2)^{3/4} (1 - \kappa \beta)^{1/2}.
\]

For the sake of completeness we state the following lemma whose
proof can be found in \cite{pon}, p. 58.

\begin{lema}
The polynomial $L(\lambda) = p_0 \lambda ^3 + p_1 \lambda ^2 + p_2
\lambda + p_3$, $p_0 > 0$, with real coefficients has all roots
with negative real parts if and only if the numbers $p_1 , p_2 ,
p_3$ are positive and the inequality $p_1 p_2 > p_0 p_3$ is
satisfied.

\label{routh}
\end{lema}

\begin{teo}
If
\begin{equation}
\varepsilon > \varepsilon_c = 2 \; \alpha \; \beta^{3/2} (1 -
\kappa \beta)^{1/2}, \label{varepsiloncritico}
\end{equation}
then the WGSS differential equations (\ref{sistemafinal}) have an
asymptotically stable equilibrium point at $P_0$. If
\[
0 < \varepsilon < \varepsilon_c
\]
then $P_0$ is unstable.

\label{teoestabilidade}
\end{teo}

\noindent {\bf Proof.} The characteristic polynomial of $Df \left(
P_0 \right)$ is given by $p(\lambda)$, where
\[
-p(\lambda)= \lambda^3 + p_1 \: \lambda^2 + p_2 \: \lambda + p_3,
\]
\[
p_1 = \varepsilon, \: \: p_2 = \frac{1 - \beta^2}{\beta}, \:\: p_3 =
\frac{2 \alpha \beta^{3/2} (1 - \beta ^2) (1 - \kappa \beta)^{1/2}
}{\beta}.
\]
The coefficients of $- p(\lambda)$ are positive. Thus a necessary
and sufficient condition for the asymptotic stability of the
equilibrium point $P_0$, as provided by the condition for one real
negative root and a pair of complex conjugate roots with negative
real part, is given by (\ref{varepsiloncritico}), according to
Lemma \ref{routh}.
\begin{flushright}
$\blacksquare$
\end{flushright}

In terms of the WGSS physical parameters, condition
(\ref{varepsiloncritico}) is equivalent to
\begin{equation}
\frac{b \: I}{m} \: \eta > 1, \label{Vichneg}
\end{equation}
where
\begin{equation}
\eta = \left| \frac{d \Omega_0}{dF} \right| = \frac{1}{2
\beta^{3/2} (1 - \kappa \beta)^{1/2}} \label{naouniformidade}
\end{equation}
is the non-uniformity of the performance of the engine which
quantifies the change in the engine speed with respect to the load
(see \cite{pon}, p. 219, for more details). Eq.
(\ref{naouniformidade}) can be easily written in terms of the
original parameters.

The rules formulated by Vyshnegradskii to enhance the stability
follow directly from (\ref{Vichneg}). In particular, the
interpretation of (\ref{Vichneg}) is that a sufficient amount of
damping $b$ must be present relative to the other physical
parameters for the system to be stable at the desired operating
speed. The condition (\ref{Vichneg}) is equivalent to the original
condition given by Vyshnegradskii for the WGS (see \cite{pon}, p.
219).

In section \ref{S4} we study the stability of $P_0$ under the
condition
\begin{equation}
\varepsilon = \varepsilon_c, \label{valorcritico}
\end{equation}
\noindent that is, on the hypersurface ---the Hopf hypersurface---
complementary to the range of validity of Theorem
\ref{teoestabilidade}.

\section{Lyapunov coefficients}\label{S3}

The beginning of this section is a review of the method found in
\cite{kuznet}, pp 177-181, and  in  \cite{kuznet2} for the
calculation of the first and second Lyapunov coefficients. The
calculation of the third Lyapunov coefficient can be found in
\cite{smb2}. The calculation of the fourth Lyapunov coefficient
has not been found by the authors in the current literature. The
extensive calculations and the  long expressions for these
coefficients have been corroborated with the software MATHEMATICA
5 \cite{math}.

Consider the differential equations
\begin{equation}
{\bf x}' = f ({\bf x}, {\bf \mu}), \label{diffequat}
\end{equation}
where ${\bf x} \in \R^n$ and ${\bf \mu} \in \R^m$ are respectively
vectors representing phase variables and control parameters.
Assume that $f$ is of class $C^{\infty}$ in $\R^n \times \R^m$.
Suppose (\ref{diffequat}) has an equilibrium point ${\bf x} = {\bf
x_0}$ at ${\bf \mu} = {\bf \mu_0}$ and, denoting the variable
${\bf x}-{\bf x_0}$ also by ${\bf x}$, write
\begin{equation}
F({\bf x}) = f ({\bf x}, {\bf \mu_0}) \label{Fhomo}
\end{equation}
as
{\small
\begin{eqnarray}\label{taylorexp}
F({\bf x}) = A{\bf x} + \frac{1}{2} \: B({\bf x},{\bf x}) +
\frac{1}{6} \: C({\bf x}, {\bf x}, {\bf x}) + \: \frac{1}{24} \:
D({\bf x}, {\bf x}, {\bf x}, {\bf x}) + \frac{1}{120} \: E({\bf
x}, {\bf x}, {\bf x}, {\bf x}, {\bf x}) + {\nonumber} \\
\frac{1}{720} \: K({\bf x}, {\bf x}, {\bf x}, {\bf x}, {\bf x},
{\bf x})+ \frac{1}{5040} \: L({\bf x}, {\bf x}, {\bf x}, {\bf x},
{\bf x}, {\bf x}, {\bf x}) +  \frac{1}{40320} \: M({\bf x}, {\bf
x}, {\bf x}, {\bf x}, {\bf x}, {\bf x}, {\bf x},{\bf x})  \\+
\frac{1}{362880} \: N({\bf x}, {\bf x}, {\bf x}, {\bf x}, {\bf x},
{\bf x}, {\bf x},{\bf x},{\bf x}) + O(|| {\bf x}
||^{10}){\nonumber},
\end{eqnarray}
}
\noindent where $A = f_{\bf x}(0,{\bf \mu_0})$ and
{\small
\begin{equation}
B_i ({\bf x},{\bf y}) = \sum_{j,k=1}^n \frac{\partial ^2
F_i(\xi)}{\partial \xi_j \: \partial \xi_k} \bigg|_{\xi=0} x_j \;
y_k, \label{Bap}
\end{equation}
}
{\small
\begin{equation}
C_i ({\bf x},{\bf y},{\bf z}) = \sum_{j,k,l=1}^n \frac{\partial ^3
F_i(\xi)}{\partial \xi_j \: \partial \xi_k \: \partial \xi_l}
\bigg|_{\xi=0} x_j \; y_k \: z_l, \label{Cap}
\end{equation}
}
{\small
\begin{equation}
D_i ({\bf x},{\bf y},{\bf z},{\bf u}) = \sum_{j,k,l,r=1}^n
\frac{\partial ^4 F_i(\xi)}{\partial \xi_j \: \partial \xi_k \:
\partial \xi_l \: \partial \xi_r} \bigg|_{\xi=0} x_j \;
y_k \: z_l \: u_r, \label{Dap}
\end{equation}
}
{\small
\begin{equation}
E_i ({\bf x},{\bf y},{\bf z},{\bf u},{\bf v}) = \sum_{j,k,l,r,p
=1}^n \frac{\partial ^5 F_i(\xi)}{\partial \xi_j \: \partial \xi_k
\: \partial \xi_l \: \partial \xi_r \: \partial \xi_p}
\bigg|_{\xi=0} x_j \; y_k \: z_l \: u_r \: v_p, \label{Eap}
\end{equation}
}
{\small
\begin{equation}
K_i ({\bf x},{\bf y},{\bf z},{\bf u},{\bf v},{\bf w}) =
\sum_{j,\ldots,q =1}^n \frac{\partial ^6 F_i(\xi)}{\partial \xi_j
\: \partial \xi_k \: \partial \xi_l \: \partial \xi_r \: \partial
\xi_p \: \partial \xi_q} \bigg|_{\xi=0} x_j \; y_k \: z_l \: u_r
\: v_p \: w_q, \label{Kap}
\end{equation}
}
{\small
\begin{equation}
L_i ({\bf x},{\bf y},{\bf z},{\bf u},{\bf v},{\bf w},{\bf t}) =
\sum_{j,\ldots,h =1}^n \frac{\partial ^7 F_i(\xi)}{\partial \xi_j
\partial \xi_k
\partial \xi_l  \partial \xi_r  \partial \xi_p  \partial \xi_q \partial \xi_h} \bigg|_{\xi=0} x_j \;
y_k \: z_l \: u_r \: v_p \: w_q \: t_h, \label{Lap}
\end{equation}
}
{\small
\begin{equation}
M_i ({\bf x},{\bf y},{\bf z},{\bf u},{\bf v},{\bf w},{\bf t},{\bf
r}) = \sum_{j,\ldots,a =1}^n \frac{\partial ^8 F_i(\xi)}{\partial
\xi_j
\ldots
\partial \xi_h \partial \xi_a} \bigg|_{\xi=0}
x_j \: y_k \: z_l \: u_r \: v_p \: w_q \: t_h \: r_a, \label{Map}
\end{equation}
}
{\small
\begin{equation}
N_i ({\bf x},{\bf y},{\bf z},{\bf u},{\bf v},{\bf w},{\bf t},{\bf
r},{\bf s}) = \sum_{j,\ldots,b =1}^n \frac{\partial ^9
F_i(\xi)}{\partial \xi_j
\ldots
\partial \xi_b} \bigg|_{\xi=0} x_j \;
y_k \: z_l \: u_r \: v_p \: w_q \: t_h \: r_a \: s_b, \label{Nap}
\end{equation}
}
\noindent for $i = 1, \ldots, n$.

Suppose $({\bf x_0}, {\bf \mu_0})$ is an equilibrium point of
(\ref{diffequat}) where the Jacobian matrix $A$ has a pair of
purely imaginary eigenvalues $\lambda_{2,3} = \pm i \omega_0$,
$\omega_0 > 0$, and admits no other eigenvalue with zero real
part. Let $T^c$ be the generalized eigenspace of $A$ corresponding
to $\lambda_{2,3}$. By this is meant that it is the largest
subspace invariant by $A$ on which the eigenvalues are
$\lambda_{2,3}$.

Let $p, q \in \mathbb C ^n$ be vectors such that
\begin{equation}
A q = i \omega_0 \: q,\:\: A^{\top} p = -i \omega_0 \: p, \:\:
\langle p,q \rangle = \sum_{i=1}^n \bar{p}_i \: q_i \:\: = 1,
\label{normalization}
\end{equation}
where $A^{\top}$ is the transposed matrix. Any vector $y \in T^c$
can be represented as $y = w q + \bar w \bar q$, where $w =
\langle p , y \rangle \in \mathbb C$. The two dimensional center
manifold can be parameterized by $w , \bar w$, by means of an
immersion of the form  ${\bf x} = H (w, \bar w)$, where $H:\mathbb
C^2 \to \R^n$ has a Taylor expansion of the form
\begin{equation}
H(w,{\bar w}) = w q + {\bar w}{\bar q} + \sum_{2 \leq j+k \leq 9}
\frac{1}{j!k!} \: h_{jk}w^j{\bar w}^k + O(|w|^{10}), \label{defH}
\end{equation}
with $h_{jk} \in \mathbb C ^n$ and  $h_{jk}={\bar h}_{kj}$.
Substituting this expression into (\ref{diffequat}) we obtain the
following differential equation
\begin{equation} \label{ku}
H_w w' + H_{\bar w} {\bar w}' = F (H(w,{\bar w})),
\end{equation}
where $F$ is given by (\ref{Fhomo}).

The complex vectors $h_{ij}$ are obtained solving the system of
linear equations defined by the coefficients of (\ref{ku}), taking
into account the coefficients of $F$, so that system (\ref{ku}),
on the chart $w$ for a central manifold, writes as follows {\small
\[
w'= i \omega_0 w + \frac{1}{2} \; G_{21} w |w|^2 + \frac{1}{12} \;
G_{32} w |w|^4 + \frac{1}{144} \; G_{43} w |w|^6 + \frac{1}{2880}
\; G_{54} w |w|^8 + O(|w|^{10}),
\]
}
with $G_{jk} \in \mathbb C$.

The {\it first Lyapunov coefficient} $l_1$ is defined by
\begin{equation}
l_1 =  \frac{1}{2} \: {\rm Re} \; G_{21}, \label{defcoef}
\end{equation}
where
\[
G_{21}= \langle p, \mathcal H_{21} \rangle, \; {\mbox {and}} \;
\mathcal H_{21} = C(q,q,\bar q) + B(\bar q, h_{20}) + 2 B(q,
h_{11}).
\]

The complex vector $h_{21}$ can be found by solving the
nonsingular $(n+1)$-dimensional system
\[
\left( \begin{array}{cc}
i \omega_0 I_n -A & q \\
\\ {\bar p} & 0
\end{array} \right) \left( \begin{array} {c}
h_{21}\\
\\s \end{array} \right)= \left( \begin{array}{c} \mathcal H_{21} -G_{21} q \\
\\0 \end{array} \right),\label{h20}
\]
with the condition $\langle p, h_{21} \rangle = 0$. See Remark 3.1
of \cite{smb3}. The procedure above can be adapted in connection
with the determination of $h_{32}$ and $h_{43}$.

Defining $\mathcal H_{32}$ as
\begin{eqnarray}\label{H32}
\mathcal H_{32} = 6 B(h_{11},h_{21})+ B({\bar h}_{20},h_{30}) + 3
B({\bar h}_{21},h_{20})+ 3 B(q,h_{22}) \nonumber \\ + 2 B(\bar q,
h_{31}) +6 C(q,h_{11},h_{11}) + 3 C(q, {\bar h}_{20}, h_{20})+ 3
C(q,q,{\bar h}_{21}) \nonumber \\ +6 C(q,\bar q, h_{21}) + 6
C(\bar q, h_{20}, h_{11}) + C(\bar q, \bar q, h_{30}) +
D(q,q,q,{\bar h}_{20}) \nonumber \\ + 6 D(q,q,\bar q,h_{11}) + 3
D(q, \bar q,\bar q, h_{20}) + E(q,q,q,\bar q,\bar q) \nonumber \\
-6 G_{21}h_{21} - 3 {\bar G}_{21} h_{21}, \nonumber
\end{eqnarray}
and from the coefficients of the terms $w^3 {\bar w}^2$ in
(\ref{ku}), one has a singular system for $h_{32}$
\[
(i \omega_0 I_n -A)h_{32}= \mathcal H_{32} - G_{32}q, \label{h32m}
\]
which has solution if and only if
\begin{equation}
\langle p, \mathcal H_{32} - G_{32}q \rangle = 0. \label{H32m}
\end{equation}

The {\it second Lyapunov coefficient} is defined by
\begin{equation}
l_2= \frac{1}{12} \: {\rm Re} \: G_{32}, \label{defcoef2}
\end{equation}
where, from (\ref{H32m}), $G_{32}=\langle p, \mathcal H_{32}
\rangle$.

The complex vector $h_{32}$ can be found solving the nonsingular
$(n+1)$-dimensional system
\[
\left( \begin{array}{cc}
i \omega_0 I_n -A & q \\
\\ {\bar p} & 0
\end{array} \right) \left( \begin{array} {c}
h_{32}\\
\\s \end{array} \right)= \left( \begin{array}{c} \mathcal H_{32} -G_{32} q \\
\\0 \end{array} \right),\label{h32}
\]
with the condition $\langle p, h_{32} \rangle = 0$.

Defining $\mathcal H_{43} $ as
{\small
\begin{eqnarray}\label{H43m}
\mathcal H_{43} = 12 B(h_{11},h_{32}) +
6 B(h_{20},{\bar h}_{32}) + 3 B({\bar h}_{20},h_{41}) \nonumber \\
+ 18 B(h_{21},h_{22}) + 12 B({\bar h}_{21},h_{31}) + 4
B(h_{30},{\bar h}_{31}) + B({\bar h}_{30},h_{40}) \nonumber \\ + 4
B(q,h_{33}) + 3 B(\bar q, h_{42}) + 36 C(h_{11},h_{11},h_{21}) +
36 C(h_{11},h_{20},{\bar h}_{21}) \nonumber \\ + 12 C(h_{11},{\bar
h}_{20},h_{30}) + 3 C(h_{20},h_{20},{\bar h}_{30}) + 18
C(h_{20},{\bar h}_{20},h_{21}) \nonumber \\ + 36
C(q,h_{11},h_{22}) + 12 C(q,h_{20},{\bar h}_{31}) + 12 C(q,{\bar
h}_{20},h_{31}) \nonumber \\ + 36 C(q,h_{21},{\bar h}_{21}) + 4
C(q,h_{30},{\bar h}_{30}) + 6 C(q,q,{\bar h}_{32}) \nonumber \\ +
12 C(q,\bar q,h_{32}) + 24 C(\bar q,h_{11},h_{31}) + 18 C(\bar
q,h_{20},h_{22}) \nonumber \\ + 3 C(\bar q,{\bar h}_{20},h_{40}) +
18 C(\bar q,h_{21},h_{21}) + 12 C(\bar q,{\bar h}_{21},h_{30})
\nonumber \\ + 3 C(\bar q, \bar q, h_{41}) + 24
D(q,h_{11},h_{11},h_{11}) + 36 D(q,h_{11},h_{20},{\bar h}_{20})
\nonumber \\ + 36 D(q,q,h_{11},{\bar h}_{21}) + 6
D(q,q,h_{20},{\bar h}_{30}) + 18 D(q,q,{\bar h}_{20},h_{21})
\nonumber \\ + 4 D(q,q,q,{\bar h}_{31}) + 18 D(q,q,\bar q,h_{22})
+ 72 D(q,\bar q, h_{11},h_{21}) \nonumber \\ + 36 D(q,\bar
q,h_{20},{\bar h}_{21}) + 12 D(q,\bar q,{\bar h}_{20},h_{30}) + 12
D(q,\bar q,\bar q, h_{31}) \nonumber \\ + 36 D(\bar
q,h_{11},h_{11},h_{20}) + 9 D(\bar q,h_{20},h_{20},{\bar h}_{20})
+ 12 D(\bar q,\bar q,h_{11},h_{30}) \nonumber \\ + 18 D(\bar
q,\bar q,h_{20},h_{21}) + D(\bar q,\bar q, \bar q,h_{40}) + 12
E(q,q,q,h_{11},{\bar h}_{20}) \nonumber \\ + E(q,q,q,q,{\bar
h}_{30}) + 12 E(q,q,q,\bar q,{\bar h}_{21}) + 36 E(q,q,\bar
q,h_{11},h_{11}) \nonumber \\ + 18 E(q,q,\bar q,h_{20},{\bar
h}_{20}) + 18 E(q,q,\bar q,\bar q, h_{21}) + 36 E(q,\bar q,\bar
q,h_{11},h_{20}) \nonumber \\ + 4 E(q,\bar q,\bar q,\bar q,
h_{30})  + 3 E(\bar q, \bar q, \bar q, h_{20},h_{20}) + 3
K(q,q,q,q,\bar q,{\bar h}_{20}) \nonumber \\ + 12 K(q,q,q,\bar
q,\bar q, h_{11}) + 6 K(q,q,\bar q,\bar q,\bar q,h_{20}) +
L(q,q,q,q,\bar q,\bar q,\bar q) \nonumber \\ - 6 (2G_{32}h_{21} +
{\bar G}_{32}h_{21} + 3 G_{21} h_{32} + 2 {\bar G}_{21}h_{32})
\nonumber,
\end{eqnarray}
}
and from the coefficients of the terms $w^4 {\bar w}^3$ in
(\ref{ku}), one has a singular system for $h_{43}$
\[
(i \omega_0 I_n -A)h_{43}= \mathcal H_{43}- G_{43} q
\]
which has solution if and only if
\begin{eqnarray}\label{h43}
\langle p, \mathcal H_{43}- G_{43} q \rangle =0.
\end{eqnarray}

The {\it third Lyapunov coefficient} is defined by
\begin{equation}
l_3= \frac{1}{144} \: {\rm Re} \: G_{43}, \label{defcoef3}
\end{equation}
where, from (\ref{h43}), $ G_{43} = \langle p, \mathcal H_{43}
\rangle$.

The complex vector $h_{43}$ can be found solving the nonsingular
$(n+1)$-dimensional system
\[
\left( \begin{array}{cc}
i \omega_0 I_n -A & q \\
\\ {\bar p} & 0
\end{array} \right) \left( \begin{array} {c}
h_{43}\\
\\s \end{array} \right)= \left( \begin{array}{c} \mathcal H_{43} -G_{43} q \\
\\0 \end{array} \right),\label{h32}
\]
with the condition $\langle p, h_{43} \rangle = 0$.

Defining $\mathcal H_{54} $ by the below expression

{\footnotesize
\begin{eqnarray*}
20 B(h_{11},h_{43}) + 10 B(h_{20},\bar h_{43}) + 6 B(\bar
h_{20},h_{52}) + 40 B(h_{21},h_{33}) +30 B(\bar h_{21},h_{42}) + \\
60 B(h_{22},h_{32}) +  10 B(h_{30},\bar h_{42}) + 4 B(\bar
h_{30},h_{51}) + 40 B(h_{31},\bar h_{32}) + 20 B(\bar
h_{31},h_{41}) + \\ 5 B(h_{40},\bar h_{41}) + B(\bar
h_{40},h_{50}) + 5 B(q,h_{44}) + 4 B(\bar q,h_{53}) + 120
C(h_{11},h_{11},h_{32}) + \\ 60 C(h_{11},\bar h_{20},h_{41}) + 360
C(h_{11},h_{21},h_{22}) + 240 C(h_{11},\bar h_{21},h_{31}) +  80
C(h_{11},h_{30},\bar h_{31}) + \\ 20 C(h_{11},\bar h_{30},h_{40})
+ 120 C(h_{20},h_{11},\bar h_{32}) + 15 C(h_{20},h_{20},\bar
h_{41}) + 60 C(h_{20},\bar h_{20},h_{32}) + \\ 120
C(h_{20},h_{21},\bar h_{31}) + 180 C(h_{20},\bar h_{21},h_{22}) +
10 C(h_{20},h_{30},\bar h_{40}) + 40 C(h_{20},\bar h_{30},h_{31})
+ \\ 3 C(\bar h_{20},\bar h_{20},h_{50}) + 120 C(\bar
h_{20},h_{21},h_{31}) + 30 C(\bar h_{20},\bar h_{21},h_{40}) + 60
C(\bar h_{20},h_{30},h_{22}) + \\ 180 C(h_{21},h_{21},\bar h_{21})
+ 60 C(\bar h_{21},\bar h_{21},h_{30}) + 40 C(h_{30},h_{21},\bar
h_{30}) + 80 C(q,h_{11},h_{33}) + \\ 30 C(q,h_{20},\bar h_{42}) +
30 C(q,\bar h_{20},h_{42}) + 120 C(q,h_{21},\bar h_{32}) + 120
C(q,\bar h_{21},h_{32}) + \\ 90 C(q,h_{22},h_{22}) + 20
C(q,h_{30},\bar h_{41}) + 20 C(q,\bar h_{30},h_{41}) + 80
C(q,h_{31},\bar h_{31}) + 5 C(q,h_{40},\bar h_{40}) + \\ 10
C(q,q,\bar h_{43}) + 20 C(q,\bar q,h_{43}) + 60 C(\bar
q,h_{11},h_{42}) + 40 C(\bar q,h_{20},h_{33}) + 12 C(\bar q,\bar
h_{20},h_{51}) + \\ 120 C(\bar q,h_{21},h_{32}) + 60 C(\bar q,\bar
h_{21},h_{41}) + 40 C(\bar q,h_{30},\bar h_{32}) + 4 C(\bar q,\bar
h_{30},h_{50}) + \\ 120 C(\bar q,h_{31},h_{22}) + 20 C(\bar
q,h_{40},\bar h_{31}) + 6 C(\bar q,\bar q,h_{52}) + 240
D(h_{11},h_{11},h_{11},h_{21}) + \\ 120 D(h_{11},h_{11},\bar
h_{20},h_{30}) + 360 D(h_{20},h_{11},h_{11},\bar h_{21}) + 360
D(h_{20},h_{11},\bar h_{20},h_{21}) + \\ 60
D(h_{20},h_{20},h_{11},\bar h_{30}) + 90 D(h_{20},h_{20},\bar
h_{20},\bar h_{21}) + 30 D(h_{20},\bar h_{20},\bar h_{20},h_{30})
+ \\ 360 D(q,h_{11},h_{11},h_{22}) + 240 D(q,h_{11},\bar
h_{20},h_{31}) + 720 D(q,h_{11},h_{21},\bar h_{21}) + \\ 80
D(q,h_{11},h_{30},\bar h_{30}) + 240 D(q,h_{20},h_{11},\bar
h_{31}) + 15 D(q,h_{20},h_{20},\bar h_{40}) + \\ 180
D(q,h_{20},\bar h_{20},h_{22}) + 120 D(q,h_{20},h_{21},\bar
h_{30}) + 180 D(q,h_{20},\bar h_{21},\bar h_{21}) + \\ 15 D(q,\bar
h_{20},\bar h_{20},h_{40}) + 180 D(q,\bar h_{20},\bar
h_{20},h_{21}) + 120 D(q,\bar h_{20},h_{30},\bar h_{21}) + \\ 120
D(q,q,h_{11},\bar h_{32}) + 30 D(q,q,h_{20},\bar h_{41}) + 60
D(q,q,\bar h_{20},h_{32}) + 120 D(q,q,h_{21},\bar h_{31}) + \\ 180
D(q,q,\bar h_{21},h_{22}) + 10 D(q,q,h_{30},\bar h_{40}) + 40
D(q,q,\bar h_{30},h_{31}) + 10 D(q,q,q,\bar h_{42}) + \\ 40
D(q,q,\bar q,h_{33}) + 240 D(q,\bar q,h_{11},h_{32}) + 120
D(q,\bar q,h_{20},\bar h_{32}) + 60 D(q,\bar q,\bar h_{20},h_{41})
+ \\ 360 D(q,\bar q,h_{21},h_{22}) + 240 D(q,\bar q,\bar
h_{21},h_{31}) + 80 D(q,\bar q,h_{30},\bar h_{31}) + 20 D(q,\bar
q,\bar h_{30},h_{40}) + \\ 30 D(q,\bar q,\bar q,h_{42}) + 240
D(\bar q,h_{11},h_{11},h_{31}) + 60 D(\bar q,h_{11},\bar
h_{20},h_{40}) + 360 D(\bar q,h_{11},h_{21},h_{21}) + \\ 240
D(\bar q,h_{11},h_{30},\bar h_{21}) + 360 D(\bar
q,h_{20},h_{11},h_{22}) + 60 D(\bar q,h_{20},h_{20},\bar h_{31}) +
\\ 120 D(\bar q,h_{20},\bar h_{20},h_{31}) + 360 D(\bar
q,h_{20},h_{21},\bar h_{21}) + 40 D(\bar q,h_{20},h_{30},\bar
h_{30}) +  \\ 120 D(\bar q,\bar h_{20},h_{30},h_{21}) + 60 D(\bar
q,\bar q,h_{11},h_{41}) + 60 D(\bar q,\bar q,h_{20},h_{32}) + 6
D(\bar q,\bar q,h_{20},h_{50}) + \\ 120 D(\bar q,\bar
q,h_{21},h_{31}) + 30 D(\bar q,\bar q,\bar h_{1},h_{40}) + 60
D(\bar q,\bar q,h_{30},h_{22}) + 4 D(\bar q,\bar q,\bar q,h_{51})
+ \\ 120 E(q,h_{11},h_{11},h_{11},h_{11}) + 360
E(q,h_{20},h_{11},h_{11},\bar h_{20}) + 45 E(q,h_{20},h_{20},\bar
h_{20},\bar h_{20}) + \\ 360 E(q,q,h_{11},h_{11},\bar h_{21}) +
360 E(q,q,h_{11},\bar h_{20},h_{21}) + 120
E(q,q,h_{20},h_{11},\bar h_{30}) + \\ 180 E(q,q,h_{20},\bar
h_{20},\bar h_{21}) + 30 E(q,q,\bar h_{20},\bar h_{20},h_{30}) +
80 E(q,q,q,h_{11},\bar h_{31}) + \\ 10 E(q,q,q,h_{20},\bar h_{40})
+ 60 E(q,q,q,\bar h_{20},h_{22}) + 40 E(q,q,q,h_{21},\bar h_{30})
+ 60 E(q,q,q,\bar h_{21},\bar h_{21}) + \\ 5 E(q,q,q,q,\bar
h_{41}) + 40 E(q,q,q,\bar q,\bar h_{32}) + 360 E(q,q,\bar
q,h_{11},h_{22}) + 120 E(q,q,\bar q,h_{20},\bar h_{31}) + \\ 120
E(q,q,\bar q,\bar h_{20},h_{31}) + 360 E(q,q,\bar q,h_{21},\bar
h_{21}) + 40 E(q,q,\bar q,h_{30},\bar h_{30}) + 60 E(q,q,\bar
q,\bar q,h_{32}) + \\ 720 E(q,\bar q,h_{11},h_{11},h_{21}) + 240
E(q,\bar q,h_{11},\bar h_{20},h_{30}) + 720 E(q,\bar
q,h_{20},h_{11},\bar h_{21}) + \\ 60 E(q,\bar q,h_{20},h_{20},\bar
h_{30}) + 360 E(q,\bar q,h_{20},\bar h_{20},h_{21}) + 240 E(q,\bar
q,\bar q,h_{11},h_{31}) + \\ 180 E(q,\bar q,\bar q,h_{20},h_{22})
+ 30 E(q,\bar q,\bar q,\bar h_{20},h_{40}) + 180 E(q,\bar q,\bar
q,h_{21},h_{21}) + 120 E(q,\bar q,\bar q,h_{30},\bar h_{21}) + \\
20 E(q,\bar q,\bar q,\bar q,h_{41}) + 240 E(\bar
q,h_{20},h_{11},h_{11},h_{11}) + 180 E(\bar
q,h_{20},h_{20},h_{11},\bar h_{20}) + \\ 120 E(\bar q,\bar
q,h_{11},h_{11},h_{30}) + 360 E(\bar q,\bar
q,h_{20},h_{11},h_{21}) + 90 E(\bar q,\bar q,h_{20},h_{20},\bar
h_{21}) + \\ 60 E(\bar q,\bar q,h_{20},\bar h_{20},h_{30}) + 20
E(\bar q,\bar q,\bar q,h_{11},h_{40}) + 40 E(\bar q,\bar q,\bar
q,h_{20},h_{31}) + 40 E(\bar q,\bar q,\bar q,h_{30},h_{21}) + \\
E(\bar q,\bar q,\bar q,\bar q,h_{50}) + 120
K(q,q,q,h_{11},h_{11},\bar h_{20}) + 30 K(q,q,q,h_{20},\bar
h_{20},\bar h_{20}) + \\ 20 K(q,q,q,q,h_{11},\bar h_{30}) + 30
K(q,q,q,q,\bar h_{20},\bar h_{21}) + K(q,q,q,q,q,\bar h_{40}) + 20
K(q,q,q,q,\bar q,\bar h_{31}) + \\ 240 K(q,q,q,\bar q,h_{11},\bar
h_{21}) + 40 K(q,q,q,\bar q,h_{20},\bar h_{30}) + 120 K(q,q,q,\bar
q,\bar h_{20},h_{21}) + \\ 60 K(q,q,q,\bar q,\bar q,h_{22}) + 240
K(q,q,\bar q,h_{11},h_{11},h_{11}) + 360 K(q,q,\bar
q,h_{20},h_{11},\bar h_{20}) + \\ 360 K(q,q,\bar q,\bar
q,h_{11},h_{21}) + 180 K(q,q,\bar q,\bar q,h_{20},\bar h_{21}) +
60 K(q,q,\bar q,\bar q,\bar h_{20},h_{30}) + \\ 40 K(q,q,\bar
q,\bar q,\bar q,h_{31}) + 360 K(q,\bar q,\bar
q,h_{20},h_{11},h_{11}) + 90 K(q,\bar q,\bar q,h_{20},h_{20},\bar
h_{20}) + \\ 80 K(q,\bar q,\bar q,\bar q,h_{11},h_{30}) + 120
K(q,\bar q,\bar q,\bar q,h_{20},h_{21}) + 5 K(q,\bar q,\bar q,\bar
q,\bar q,h_{40}) + \\ 60 K(\bar q,\bar q,\bar
q,h_{20},h_{20},h_{11}) + 10 K(\bar q,\bar q,\bar q,\bar
q,h_{20},h_{30}) + 3 L(q,q,q,q,q,\bar h_{20},\bar h_{20}) + \\ 4
L(q,q,q,q,q,\bar q,\bar h_{30}) + 60 L(q,q,q,q,\bar q,h_{11},\bar
h_{20}) + 30 L(q,q,q,q,\bar q,\bar q,\bar h_{21}) + \\ 120
L(q,q,q,\bar q,\bar q,h_{11},h_{11}) + 60 L(q,q,q,\bar q,\bar
q,h_{20},\bar h_{20}) + 40 L(q,q,q,\bar q,\bar q,\bar q,h_{21}) +
\\ 120 L(q,q,\bar q,\bar q,\bar q,h_{20},h_{11}) + 10 L(q,q,\bar q,\bar q,\bar q,\bar
q,h_{30}) + 15 L(q,\bar q,\bar q,\bar q,\bar q,h_{20},h_{20}) + \\
6 M(q,q,q,q,q,\bar q,\bar q,\bar h_{20}) + 20 M(q,q,q,q,\bar
q,\bar q,\bar q,h_{11}) + 10 M(q,q,q,\bar q,\bar q,\bar q,\bar
q,h_{20}) + \\ N(q,q,q,q,q,\bar q,\bar q,\bar q,\bar q),
\end{eqnarray*}
}
and from the coefficients of the terms $w^5 {\bar w}^4$ in
(\ref{ku}), one has a singular system for $h_{54}$
\[
(i \omega_0 I_n -A)h_{54}= \mathcal H_{54}- G_{54} q
\]
which has solution if and only if
\begin{eqnarray}\label{h54}
\langle p, \mathcal H_{54}- G_{54} q \rangle =0.
\end{eqnarray}

The {\it fourth Lyapunov coefficient} is defined by
\begin{equation}
l_4= \frac{1}{2880} \: {\rm Re} \: G_{54}, \label{defcoef4}
\end{equation}
where, from (\ref{h54}), $ G_{54} = \langle p, \mathcal H_{54}
\rangle$.

\begin{remark}\label{conceitual}
Other equivalent definitions and algorithmic  procedures to write
the expressions for the Lyapunov coefficients $l_j , j= 1,2,3,4$,
for two dimensional systems can be found in Andronov et al.
\cite{al} and Gasull et al. \cite{gt}, among others. These
procedures apply also to the three dimensional systems of this
work, if  properly restricted to the center manifold. The authors
found, however, that the  method  outlined above, due to Kuznetsov
\cite{kuznet, kuznet2}, requiring  no explicit formal evaluation
of the center manifold, is better adapted to the needs of this
work.
\end{remark}

A {\it Hopf point} $({\bf x_0}, {\bf \mu_0})$ is an equilibrium
point of (\ref{diffequat}) where the Jacobian matrix $A = f_{\bf
x}({\bf x_0}, {\bf \mu_0})$ has a pair of purely imaginary
eigenvalues $\lambda_{2,3} = \pm i \omega_0$, $\omega_0 > 0$, and
admits  no other critical eigenvalues ---i.e. located on the
imaginary axis. At a Hopf point a two dimensional center manifold
is well-defined, it is invariant under the flow generated by
(\ref{diffequat}) and can be continued with arbitrary high class
of differentiability to nearby parameter values. In fact, what is
well defined is the $\infty$-jet ---or infinite Taylor series---
of the center manifold, as well as that of its continuation, any
two of them having contact in the arbitrary high  order of their
differentiability class.

A Hopf point is called {\it transversal} if the parameter
dependent complex eigenvalues cross the imaginary axis with
non-zero derivative. In a neighborhood of a transversal Hopf point
---H1 point, for concision--- with $l_1 \neq 0$ the dynamic
behavior of the system (\ref{diffequat}), reduced to the family of
parameter-dependent continuations of the center manifold, is
orbitally topologically equivalent to the following complex normal
form
\[
w' = (\eta + i \omega) w + l_1 w |w|^2 ,
\]
$w \in \mathbb C $, $\eta$, $\omega$ and $l_1$ are real functions
having  derivatives of arbitrary  high order, which are
continuations  of $0$, $\omega_0$ and the first Lyapunov
coefficient at the H1 point. See  \cite{kuznet}. As $l_1 < 0$
($l_1 > 0$) one family of stable (unstable) periodic orbits can be
found on this family of manifolds, shrinking  to an equilibrium
point at the H1 point.

A {\it Hopf point of codimension 2} is a Hopf point where $l_1$
vanishes. It is called {\it transversal} if $\eta = 0$ and $l_1 =
0$ have transversal intersections, where $\eta = \eta (\mu)$ is
the real part of the critical eigenvalues. In a neighborhood of a
transversal Hopf point of codimension 2 ---H2 point, for
concision---  with $l_2 \neq 0$ the dynamic behavior of the system
(\ref{diffequat}), reduced to the family of parameter-dependent
continuations of the center manifold, is orbitally topologically
equivalent to
\[
w' = (\eta + i \omega_0) w + \tau w |w|^2 + l_2 w |w|^4 ,
\]
where $\eta$ and $\tau$ are unfolding parameters.  See
\cite{kuznet}. The bifurcation diagrams for $l_2 \neq 0$ can be
found in \cite{kuznet}, p. 313, and in \cite{takens}.

A {\it Hopf point of codimension 3} is a Hopf point of codimension
2 where $l_2$ vanishes. A Hopf point of codimension 3 point is
called {\it transversal} if $\eta = 0$, $l_1 = 0$ and $l_2 = 0$
have transversal intersections. In a neighborhood of a transversal
Hopf point of codimension 3 ---H3 point, for concision--- with
$l_3 \neq 0$ the dynamic behavior of the system (\ref{diffequat}),
reduced to the family of parameter-dependent continuations of the
center manifold, is orbitally topologically equivalent to
\[
w' = (\eta + i \omega_0) w + \tau w |w|^2 + \nu w |w|^4 + l_3 w
|w|^6 ,
\]
where $\eta$, $\tau$ and $\nu$ are unfolding parameters. The
bifurcation diagram for $l_3 \neq 0$ can be found in Takens
\cite{takens} and in \cite{smb3}.

A {\it Hopf point of codimension 4} is a Hopf point of codimension
3 where $l_3$ vanishes. A Hopf point of codimension 4 is called
{\it transversal} if $\eta = 0$, $l_1 = 0$, $l_2 = 0$ and $l_3 =
0$ have transversal intersections. In a neighborhood of a
transversal Hopf point of codimension 4 ---H4 point, for
concision--- with $l_4 \neq 0$ the dynamic behavior of the system
(\ref{diffequat}), reduced to the family of parameter-dependent
continuations of the center manifold, is orbitally topologically
equivalent to
\[
w' = (\eta + i \omega_0) w + \tau w |w|^2 + \nu w |w|^4 + \sigma w
|w|^6 + l_4 w |w|^8 ,
\]
where $\eta$, $\tau$, $\nu$ and $\sigma$ are unfolding parameters.

\begin{teo}
Suppose that the system
\[
{\bf x}' = f({\bf x},{\bf \mu}), \: {\bf x}=(x,y,z), \: \mu =
(\beta, \alpha, \kappa, \varepsilon)
\]
has the equilibrium ${\bf x} = {\bf 0}$ for $\mu = 0$ with
eigenvalues
\[
\lambda_{2,3} (\mu) = \eta (\mu) \pm i \omega(\mu),
\]
where $\omega(0) = \omega_0 > 0$. For $\mu = 0$ the following
conditions hold
\[
\eta (0) = 0, \: l_1(0) = 0, \: l_2(0) = 0, \: l_3(0) = 0,
\]
where $l_1(\mu)$, $l_2(\mu)$ and $l_3(\mu)$ are the first, second
and third Lyapunov coefficients, respectively. Assume that the
following genericity conditions are satisfied
\begin{enumerate}
\item $l_4 (0) \neq 0$, where $l_4 (0)$ is the fourth Lyapunov
coefficient;

\item the map $\mu \to (\eta (\mu), l_1 (\mu), l_2 (\mu), l_3
(\mu))$ is regular at $\mu = 0$.

\end{enumerate}
Then, by the introduction of a complex variable, the above system
reduced to the family of parameter-dependent continuations of the
center manifold, is orbitally topologically equivalent to
\[
w' = (\eta + i \omega_0) w + \tau w |w|^2 + \nu w |w|^4 + \sigma w
|w|^6 + l_4 w |w|^8
\]
where $\eta$, $\tau$, $\nu$ and $\sigma$ are unfolding parameters.

\label{teoremaHopf}
\end{teo}

\section{Hopf bifurcations in the WGSS}\label{S4}

The following theorem was proved by the authors in \cite{smb3}.

\begin{figure}[!h]
\centerline{
\includegraphics[width=10cm]{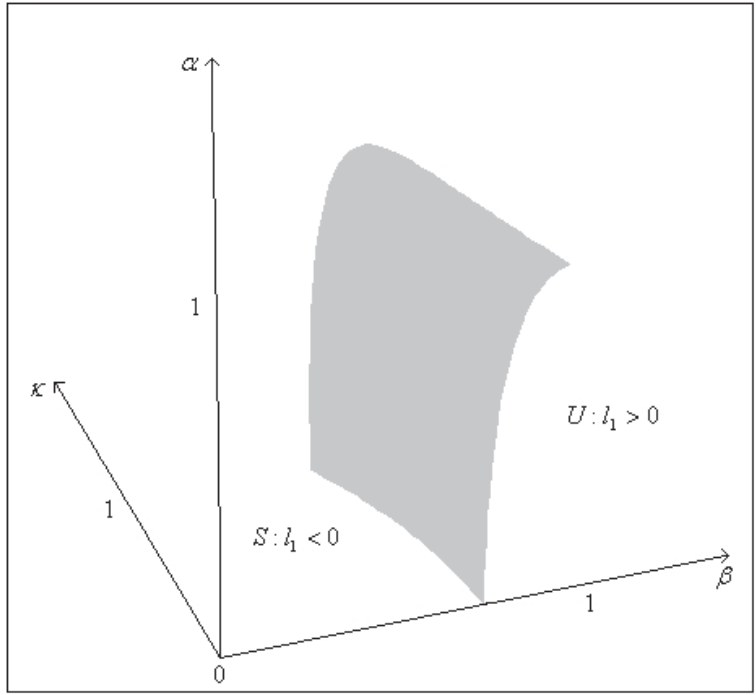}}
\caption{{\small Signs of the first Lyapunov coefficient for system
(\ref{sistemafinal})}.}

\label{sinalL1}
\end{figure}

\begin{teo}
Consider the four-parameter family of differential equations
(\ref{sistemafinal}). The first Lyapunov coefficient at the point
(\ref{P0}) for parameter values satisfying (\ref{valorcritico}) is
given by
\begin{equation}
l_1 (\beta, \alpha, \kappa) = - \frac{G_1(\beta, \alpha,
\kappa)}{4 \beta \varepsilon_c \omega_0 ^4 \omega_1 ^2
(\varepsilon_c ^4 + 5 \varepsilon_c ^2 \omega_0 ^2 + 4 \omega_0
^4)} ,\label{coeficienteL1}
\end{equation}
where
\begin{eqnarray}\label{G1}
G_1(\beta, \alpha, \kappa) = -3 + 5 \kappa \beta - (\alpha^2 -5)
\beta^2 + \kappa (\alpha ^2 -7) \beta^3 - \nonumber \\  2 \alpha^2
\kappa^2 \beta^4 - (\alpha^4 - 2 \alpha^2 \kappa^2) \beta ^6  +
\alpha ^4 \kappa \beta^7 .
\end{eqnarray}
If $G_1$ is different from zero then the system (\ref{sistemafinal})
has a transversal Hopf point at $P_0$ for $\varepsilon =
\varepsilon_c$. More specifically, if $(\beta, \alpha, \kappa) \in S
\cup U$ and $\varepsilon = \varepsilon_c$ then the system
(\ref{sistemafinal}) has an H1 point at $P_0$; if $(\beta, \alpha,
\kappa) \in S$ and $\varepsilon = \varepsilon_c$ then the H1 point
at $P_0$ is asymptotically stable and for each $\varepsilon <
\varepsilon_c$, but close to $\varepsilon_c$, there exists a stable
periodic orbit near the unstable equilibrium point $P_0$; if
$(\beta, \alpha, \kappa) \in U$ and $\varepsilon = \varepsilon_c$
then the H1 point at $P_0$ is unstable and for each $\varepsilon >
\varepsilon_c$, but close to $\varepsilon_c$, there exists an
unstable periodic orbit near the asymptotically stable equilibrium
point $P_0$. See Fig \ref{sinalL1}.

\label{teoremaHexa}
\end{teo}

\begin{figure}[!h]
\centerline{
\includegraphics[width=10cm]{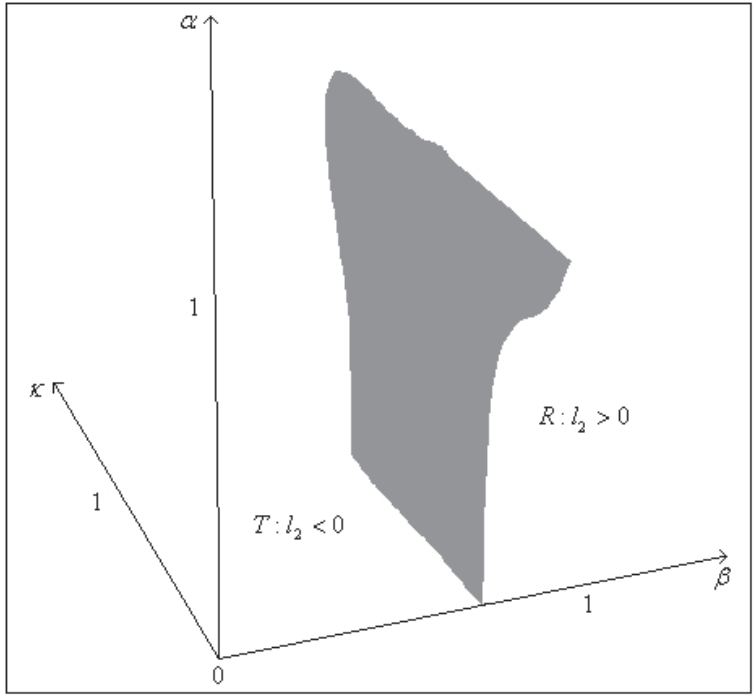}}
\caption{{\small Signs of the second Lyapunov coefficient for system
(\ref{sistemafinal})}.}

\label{sinalL2}
\end{figure}

\begin{figure}[!h]
\centerline{
\includegraphics[width=10cm]{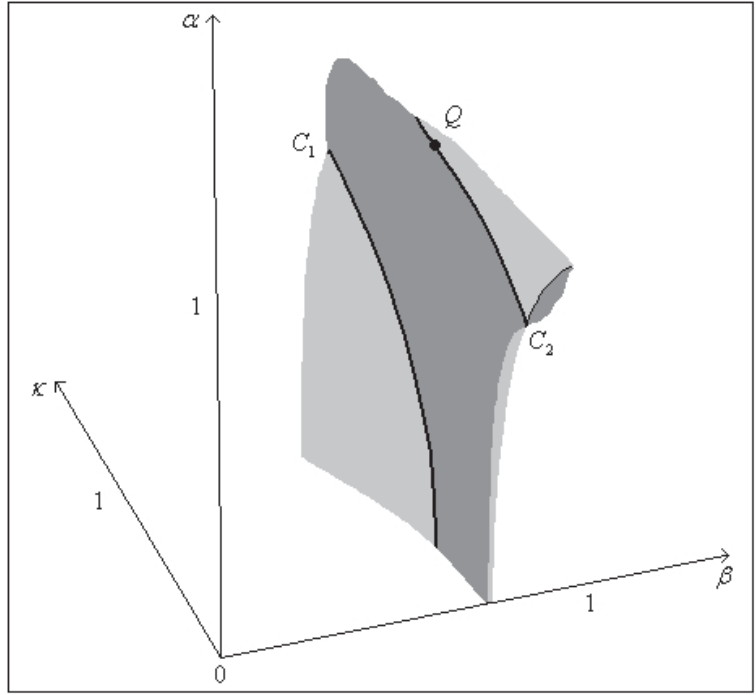}}
\caption{{\small Surfaces $l_1 = 0$ and $l_2 = 0$ and the
intersection curves}.}

\label{L1=L2=0}
\end{figure}

\begin{figure}[!h]
\centerline{
\includegraphics[width=12cm]{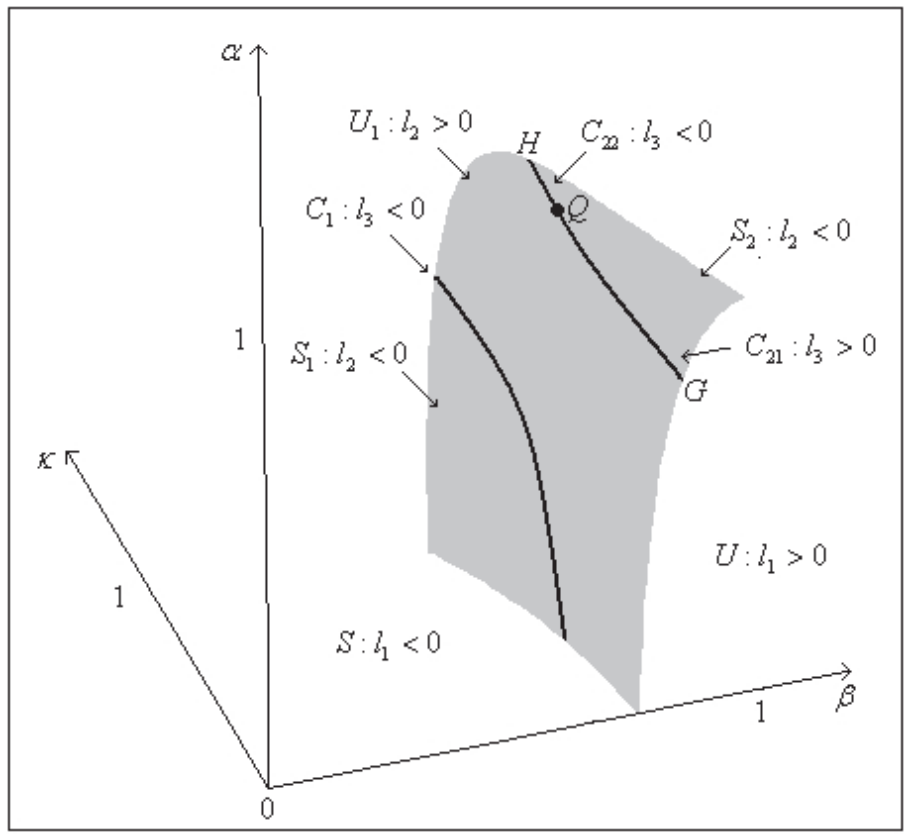}}
\caption{{\small Signs of $l_1$, $l_2$ and $l_3$}.}

\label{PointQ}
\end{figure}

\begin{figure}[!h]
\centerline{
\includegraphics[width=12cm]{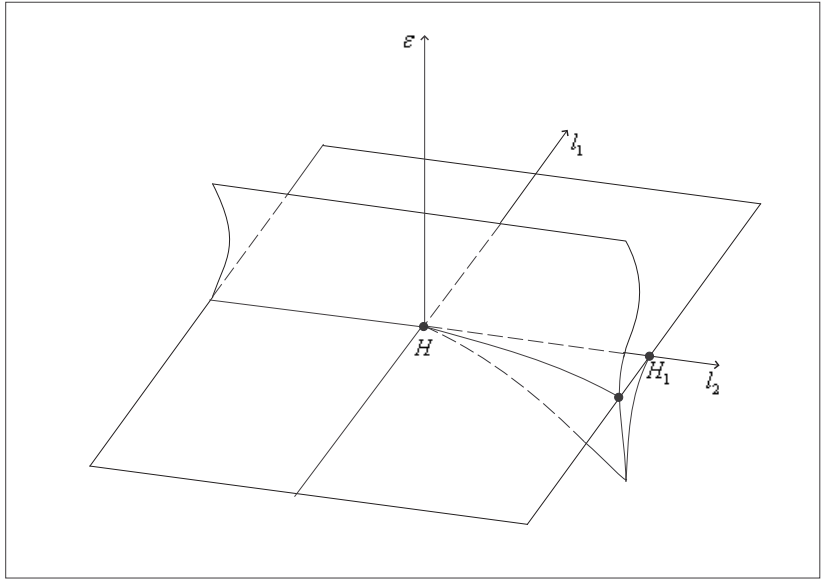}}
\caption{{\small Bifurcation diagram for a typical point $H$ where
$l_3 (H) < 0$}.}

\label{PointH}
\end{figure}

\begin{teo}
For the four-parameter family of differential equations
(\ref{sistemafinal}) there is unique point $Q = (\beta, \alpha,
\kappa, \varepsilon_c)$, with coordinates
\[
\beta = 0.93593 \ldots,\, \alpha = 1.02753 \ldots, \kappa =
0.90164 \ldots, \; \varepsilon_c = 0.73522 \ldots,
\]
where the surfaces $l_1 = 0$, $l_2 = 0 $ and $l_3 = 0 $ on the
critical hypersurface intersect and there do it transversally.
Moreover, the codimension 4 Hopf point at $P_0$ is asymptotically
stable since $l_4(Q) < 0$. More specifically, if $(\beta, \alpha,
\kappa) \in S_1 \cup S_2 \cup U_1$ and $\varepsilon = \varepsilon_c$
then the system (\ref{sistemafinal}) has an H2 point at $P_0$; if
$(\beta, \alpha, \kappa) \in S_1 \cup S_2$ and $\varepsilon =
\varepsilon_c$ then the H2 point at $P_0$ is asymptotically stable;
if $(\beta, \alpha, \kappa) \in U_1$ and $\varepsilon =
\varepsilon_c$ then the H2 point at $P_0$ is unstable. Along the
curves $C_1$ and $C_2 = C_{21} \cup C_{22} \cup \{Q\}$ of Fig
\ref{L1=L2=0} $l_2$ vanishes.If $(\beta, \alpha, \kappa) \in C_1
\cup C_{21} \cup C_{22}$ (see Fig \ref{PointQ}) and $\varepsilon =
\varepsilon_c$ then the four-parameter family of differential
equations (\ref{sistemafinal}) has a transversal Hopf point of
codimension 3 at $P_0$; if $(\beta, \alpha, \kappa) \in C_1 \cup
C_{22}$ and $\varepsilon = \varepsilon_c$ then the H3 point at $P_0$
is asymptotically stable and the bifurcation diagram for a typical
point $H$ is draw in Fig \ref{PointH}; if $(\beta, \alpha, \kappa)
\in C_{21}$ and $\varepsilon = \varepsilon_c$ then the H2 point at
$P_0$ is unstable and the bifurcation diagram for a typical point
$G$ can be found in \cite{smb2}.

\label{crucial}
\end{teo}

\noindent{\bf Computer assisted Proof.} The algebraic expression for
the second Lyapunov coefficient can be obtained in \cite{mello}.
This is too long to be put in print. The surface where the second
Lyapunov coefficient vanishes is illustrated in Fig. \ref{sinalL2}.

\begin{figure}[h!]
\hspace{-2.5cm}\centerline{
\includegraphics[width=14cm]{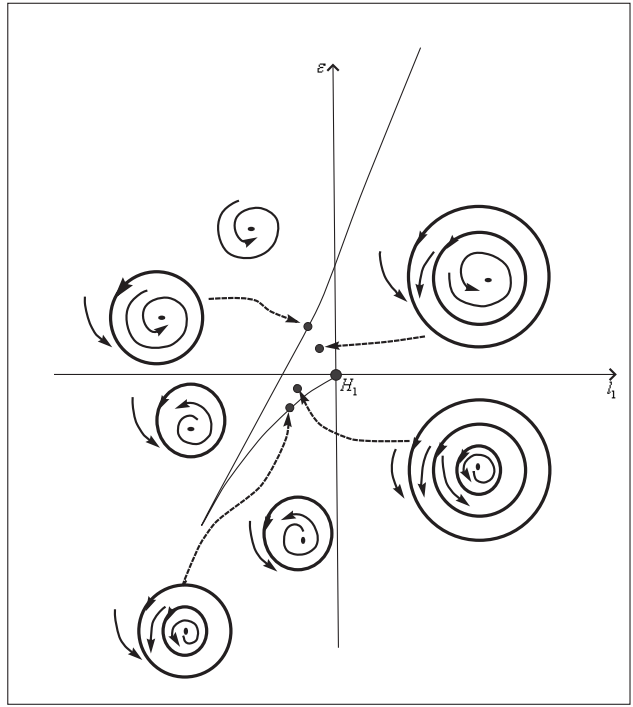}}
\caption{{\small Bifurcation diagram for a typical point $H_1$. See Fig. \ref{PointH}}.}

\label{PointH1}
\end{figure}

The intersections of the surfaces $l_1=0$ and $l_2=0$ determine the
curves $C_1$ and $C_2$ (see Fig \ref{L1=L2=0}). The signs of the
second Lyapunov coefficient on the surface $l_1 = 0$ complementary
to the curves $C_1$ and $C_2$, that is on $S_1 \cup S_2 \cup U_1$
(see Fig. \ref{PointQ}), are the following: $l_2$ is negative on
$S_1 \cup S_2$ and is positive on $U_1$ and they can be viewed as
extensions of the signs of the second Lyapunov coefficient at points
on the curve determined by the intersection of the surface $l_1 = 0$
and the plane $\kappa = 0$ studied by the authors in \cite{smb2}.
The bifurcation diagram for a typical point $G$ where $l_3 (G) > 0$
can be viewed in \cite{smb2}. In Fig \ref{PointH} and \ref{PointH1}
are illustrated the bifurcation diagrams for a typical point $H$
where $l_3 (H) < 0$.

The point $Q$ is the intersection of the surfaces $l_1 = 0 $, $l_2 =
0$ and $l_3 = 0$. The existence and uniqueness of $Q$ with the above
coordinates has been established numerically with the software
MATHEMATICA 5.

For the  point $Q$ take five decimal round-off coordinates $\beta =
0.93593$, $\alpha = 1.02753$, $\kappa = 0.90164$ and $\varepsilon_c
= 0.73522$. For these values of the parameters one has
\[
p = \left( -i/2, 0.27041 - 0.54618 i, 0.40395 + 0.20000 i \right),
\]
\[
q = \left( -i, 0.36401, 0.99407 \right),
\]
\[
h_{11}= \left( -2.65769, 0, 0.19650 \right),
\]
\[
h_{20}= \left( -4.11029 - 0.18429 i, 0.13416 -
  2.99241 i, 0.09159 - 3.36395 i \right),
\]
\[
h_{30}= \left( -3.63589 + 23.03616 i, -25.15645 -
  3.97054 i, -18.16113 - 1.69167 i \right),
\]
\begin{equation}
G_{21}= - 3.91814 i, \label{G21l3}
\end{equation}
\[
h_{21}= \left( 3.24775 + 1.67247 i, -4.52694 + 1.18222 i, 4.85950
+ 3.71541 i \right),
\]
\[
h_{40}= \left( \begin{array}{c} 160.39204 + 51.10539 i \\
    -74.41230 + 233.53975 i \\ -25.03366 + 127.34049 i
    \end{array} \right),
\]
\[
h_{31}= \left( -69.44664 - 38.56274 i, 25.90851 - 2.24484 i,
36.10391 - 65.85524 i \right),
\]
\[
h_{22} = \left( -64.50829, 0, 10.76131 \right),
\]
\begin{equation}
G_{32} = -153.21726 i, \label{G32l3}
\end{equation}
\[
h_{50}= \left( \begin{array}{c} 702.48693 - 1263.93346 i \\
    2300.44688 + 1278.57511 i \\ 1054.20770 + 363.36145 i
    \end{array} \right),
\]
\[
h_{32}= \left( \begin{array}{c} 178.24934 + 273.66781 i \\
    -233.17715 + 26.70966 i \\ 395.89053 + 272.77265 i
    \end{array} \right),
\]
\[
h_{41}= \left( \begin{array}{c} -521.71430 + 1074.26121 i \\
    -631.58388 - 484.25803 i \\ -865.10385 - 413.20000 i
    \end{array} \right),
\]
\[
h_{60}= \left( \begin{array}{c}  -10130.73267 - 9995.21750 i \\
    21830.38995 - 22126.36639 i \\ 5429.65950 - 9557.27148 i
    \end{array} \right),
\]
\[
h_{51}= \left( \begin{array}{c} 14227.43860 + 8237.49829 i \\
    -9991.87299 + 14431.55078 i \\ -5753.08267 + 11280.54380 i
    \end{array} \right),
\]
\[
h_{42}= \left( \begin{array}{c} -4351.45992 - 4936.33553 i \\
    2272.08822 + 1527.90723 i \\ 4841.97866 - 5445.36779 i
    \end{array} \right),
\]
\[
h_{33} = \left( -5969.63958, 0, 1764.47230 \right),
\]
\[
h_{70}= \left( \begin{array}{c} -146941.54096 + 63522.80004 i \\
    -161862.28504 - 374421.36634 i \\ -86069.40319 - 83969.45215 i
    \end{array} \right),
\]
\[
h_{61}= \left( \begin{array}{c} 140223.18890 - 184094.16057 i \\
    260780.07852 + 213929.28545 i \\ 151116.49070 + 92225.27059 i
    \end{array} \right),
\]
\[
h_{52}= \left( \begin{array}{c} -105557.32750 + 127994.80577 i \\
    -41289.02476 - 79039.91108 i \\ -106857.14273 - 88122.45467 i
    \end{array} \right),
\]
\begin{equation}
G_{43} = - 22328.21224 i .\label{G43l3}
\end{equation}
\[
h_{43}= \left( \begin{array}{c} 26579.27090 + 62051.16515 i \\
     -36944.56779 + 2499.10743 i \\ 78144.32459 + 54070.14624 i
    \end{array} \right),
\]
\[
h_{80}= \left( \begin{array}{c} -247681.58290 + 2173895.03048 i \\
     -6330624.44741 - 721276.35507 i \\ -1324248.15135 + 594661.38331 i
    \end{array} \right),
\]
\[
h_{71}= \left( \begin{array}{c} -2230744.30930 - 2511854.85381 i \\
     4663683.99275 - 4038564.75411 i \\ 1618564.33911 - 2037646.14488 i
    \end{array} \right),
\]
\[
h_{62}= \left( \begin{array}{c} 2540059.79128 + 2277848.86298 i \\
     -2385453.21697 + 1869088.06376 i \\ -1708253.47087 + 2025268.53034 i
    \end{array} \right),
\]
\[
h_{53}= \left( \begin{array}{c} -633499.15640 - 1125590.51413 i \\
     390598.08062 + 466226.40735 i \\ 1219484.73373 - 1101283.41903 i
    \end{array} \right),
\]
\[
h_{44} = \left( -1118100.12194, 0.00138, 546721.10946 \right),
\]
\begin{equation}
G_{54} = -22071.41115 - 5991090.52119 i .\label{G54l3}
\end{equation}
From (\ref{defcoef}), (\ref{defcoef2}), (\ref{defcoef3}),
(\ref{defcoef4}), (\ref{G21l3}), (\ref{G32l3}), (\ref{G43l3}) and
(\ref{G54l3}) one has
\[
l_1 (Q) = 0, \: l_2 (Q) = 0, \: l_3 (Q) = 0, \: l_4 (Q) =
\frac{1}{2880} \: {\rm Re} \: G_{54} = -7.66368.
\]

The calculations above have also been corroborated with 100 decimals
round-off precision performed using the software MATHEMATICA 5
\cite{math}. See \cite{mello}.

Some values of $(\alpha,\beta,\kappa) \in C_1 \cup C_2$ as well as
the corresponding values of $l_3 (\alpha,\beta,\kappa)$ are listed
in the tables below. The calculations leading to these values can be
found in \cite{mello}.

\begin{center}
\begin{tabular}{|c|c|c|c|}    \hline \hline
$\kappa$       & $\alpha$     & $\beta$  &$l_3 (\alpha,\beta,\kappa)$ on $C_1$ \\
\hline
 0.45   &  0.33319  & 0.72216  & -0.91310 \\ \hline
 0.5    &  0.42968  & 0.71770  & -0.92567 \\ \hline
 0.55   &  0.50934  & 0.71257  & -0.88152 \\ \hline
 0.6    &  0.57913  & 0.70665  & -0.82064 \\ \hline
 0.65   &  0.64241  & 0.69983  & -0.75810 \\ \hline
 0.7    &  0.70113  & 0.69201  & -0.70006 \\ \hline
 0.75   &  0.75659  & 0.68309  & -0.64900 \\ \hline
 0.8    &  0.80972  & 0.67302  & -0.60580 \\ \hline
 0.85   &  0.86120  & 0.66177  & -0.57054 \\ \hline
 0.9    &  0.91154  & 0.64940  & -0.54288 \\ \hline
 0.95   &  0.96114  & 0.63600  & -0.52217 \\ \hline\hline
\end{tabular}
\end{center}

\begin{center}
\begin{tabular}{|c|c|c|c|}    \hline \hline
$\kappa$       & $\alpha$     & $\beta$  &$l_3 (\alpha,\beta,\kappa)$ on $C_2$ \\
\hline
 0     & 0.85050   & 0.86828  & 0.39050 \\ \hline
 0.2   & 0.90524   & 0.87760  & 0.46294 \\ \hline
 0.3   & 0.93123   & 0.88397  & 0.50684 \\ \hline
 0.4   & 0.95511   & 0.89159  & 0.55538 \\ \hline
 0.5   & 0.97602   & 0.90042  & 0.60637 \\ \hline
 0.6   & 0.99330   & 0.91029  & 0.65253 \\ \hline
 0.7   & 1.00674   & 0.92071  & 0.66963 \\ \hline
 0.8   & 1.01697   & 0.93045  & 0.56860 \\ \hline
 0.9   & 1.02731   & 0.93592  & 0.01665 \\ \hline
 0.92  & 1.03020   & 0.93585  & -0.20674 \\ \hline
 0.98  & 1.04319   & 0.93201  & -1.09289 \\ \hline\hline
\end{tabular}
\end{center}

The gradients of the functions $l_1 $, $l_2 $ and $l_3$, given in
(\ref{defcoef}), (\ref{defcoef2}), (\ref{defcoef3}) at the point $Q$
are, respectively
\[
(-0.46264, 0.13437, -0.97565), (-12.44701, 2.66791, -19.19345),
\]
\[
(-266.77145, 41.80505, -372.84969).
\]
The transversality condition at $Q$ is equivalent to the
non-vanishing of the determinant of the matrix whose columns are the
above gradient vectors, which is evaluated gives $-33.31133$.

The main steps of the calculations  that provide the numerical
evidence for this theorem  have been posted in \cite{mello}.
\begin{flushright}
$\blacksquare$
\end{flushright}

\section{Concluding comments}\label{conclusion}

This paper starts reviewing the stability analysis which accounts
for the characterization, in the space of parameters, of the
structural as well as Lyapunov stability of the equilibrium of the
Watt Governor System with a Spring, WGSS. It continues with
recounting the extension of the analysis to the first order,
codimension one stable points, happening on the complement of a
surface in the critical hypersurface where the eigenvalue criterium
of Lyapunov holds, as studied by the authors \cite{smb3}, based on
the calculation of the first Lyapunov coefficient. Here the
bifurcation analysis at the equilibrium point of the WGSS is pushed
forward to the calculation of the second, third and fourth Lyapunov
coefficients which make possible the determination of the Lyapunov
as well as higher order structural stability at the equilibrium
point. See also \cite{kuznet, kuznet2}, \cite{gt} and \cite{al} .

The calculations of these coefficients, being extensive, rely  on
Computer Algebra and Numerical evaluations carried out with the
software MATHEMATICA 5 \cite{math}. In the site \cite{mello} have
been posted the main steps of the calculations in the form of
notebooks for MATHEMATICA 5.

With the analytic and numeric data provided in the analysis
performed here, the bifurcation diagrams are established along the
points of the surface where the first Lyapunov coefficient vanishes.
Pictures \ref{PointH} and \ref{PointH1} provide a qualitative synthesis of the
dynamical conclusions achieved here at the parameter values where
the WGSS achieves most complex equilibrium point. A reformulation of
these conclusions follow:

\noindent There is a ``solid tongue"  where three stable regimes
coexist: one is an equilibrium and the other two are small amplitude
periodic orbits, i.e., oscillations.

For parameters inside the ``tongue", this conclusion suggests, a
{\it hysteresis} explanation for the phenomenon of ``hunting"
observed in the performance of WGSS in an early stage of the
research on its stability conditions. Which attractor represents the
actual state of the system will depend on the path along which the
parameters evolve to reach their actual values of the parameters
under consideration. See Denny \cite{denny} for historical comments,
where he refers to the term ``hunting" to mean an oscillation around
an equilibrium going near but not reaching it.

Finally, we would like to stress that although this work ultimately
focuses the specific three dimensional, four parameter system of
differential equations given by (\ref{sistemafinal}), the method of
analysis and calculations explained in Section \ref{S3} can be
adapted to the study of other systems with three or more phase
variables and depending on four or more parameters.

\vspace{0.2cm} \noindent {\bf Acknowledgement}: The first and second
authors developed this work under the project CNPq Grants
473824/04-3 and 473747/2006-5. The first author is fellow of CNPq
and takes part in the project CNPq PADCT 620029/2004-8. The third
author is supported by CAPES. This work was finished while the
second author visited Universitat Aut\`onoma de Barcelona, supported
by CNPq grant 210056/2006-1.

\end{document}